\documentclass[runningheads]{llncs}
\usepackage[T1]{fontenc}
\usepackage{graphicx}
\usepackage{amssymb,amsmath}
\usepackage{multirow}
\usepackage{tikz}
\newcommand*\circled[2][1.6]{\tikz[baseline=(char.base)]{
    \node[shape=circle, draw, inner sep=1pt, 
        minimum height=#1em] (char) {\vphantom{WAH1g}#2};}}
\newcommand*\insquare[2][1.6]{\tikz[baseline=(char.base)]{
    \node[shape=rectangle, draw, inner sep=1pt, 
        minimum height=#1em] (char) {\vphantom{WAH1g}#2};}}
\usetikzlibrary{cd}
\usepackage[percent]{overpic}
\begin{document}
\title{An Algebraic-Symmetric Analysis of Counterpoint and Modulation in the Music of Claudio Monteverdi}
\titlerunning{An Algebraic-Symmetric Analysis...}
%
\author{Octavio A. Agustín-Aquino\inst{1}\orcidID{0000-0002-0556-6236} \and
Brandon J. Curiel-López\inst{1}}
\authorrunning{O. A. Agustín-Aquino and B. J. Curiel-López}
%
\institute{Universidad Tecnológica de la Mixteca, México
\email{octavioalberto@mixteco.utm.mx, culb990302@gs.utm.mx
}}
\maketitle              
\begin{abstract}
We address certain structural innovations in the music of Claudio Monteverdi, which defined the pivotal transition from the Renaissance \textit{prima pratica} to the Baroque \textit{seconda pratica}. To formalize this analysis, we employ Mazzola's symmetry-based framework for counterpoint and quantum/duality models for tonal modulation. Our findings demonstrate that these mathematical structures provide a rigorous validation of Monteverdi’s compositional choices, revealing an underlying logic to harmonic and contrapuntal treatments that were heavily criticized by contemporaries such as Giovanni Artusi. By quantifying concepts like compositional parsimony and modulations, we aim to provide a precise and analytical lens to the interdisciplinary field of mathematical musicology, bridging the gap between historical interpretation and formal theory.

\end{abstract}
\section{Introduction: Monteverdi and the Mathematical Turn in Music Theory}
The work of Claudio Monteverdi (1567-1643) stands at the crossroads of Western music history, marking the transition from the Renaissance to the Baroque period. This era was defined by a fervent stylistic debate between the \textit{prima pratica}, the established polyphonic style governed by strict contrapuntal rules, and the emergent \textit{seconda pratica}, a new style that prioritized textual expression over certain formal conventions. Here we show that the apparatus of Mazzola's theory for counterpoint and harmonic modulation and the $PLR/TI$ duality can explain certain features of Monteverdi's contrapuntal and harmonic language.

\section{The Musicological Context: Monteverdi at the Dawn of the Baroque}

To fully appreciate the significance of Monteverdi's innovations and the mathematical models used to describe them, it is essential to understand the historical and theoretical landscape of the late 16th and early 17th centuries. This period was characterized by a rich tradition of polyphonic music and a growing tension between established compositional rules and the demand for greater expressive power.

Monteverdi was a master of the principal genres of his time: sacred works such as the Mass and motet, the secular madrigal, and the nascent genre of opera. This versatility is significant: his controversial techniques were not mere experiments but deliberate choices of a composer fully conversant with traditional practice. It was precisely this command of the \textit{prima pratica} that gave weight to his departures from it. Our analysis reflects this breadth: we examine counterpoint in his madrigals, motets, and Mass settings, while our modulation analysis focuses on \textit{L'Orfeo (1607)} \cite{cM09}, where the demands of dramatic narrative pushed his harmonic language furthest from convention and drew the sharpest criticism.

\subsection{The \textit{prima pratica} vs. \textit{seconda pratica} Debate}

The fundamental conflict of Monteverdi's era was the stylistic schism between the \textit{prima pratica} (first practice) and the \textit{seconda pratica} (second practice).
\begin{enumerate}
\item The \textit{prima pratica}, or \textit{stilo antico}, was the established style of Renaissance polyphony. Defended by Monteverdi's brother, Giulio Cesare, it was defined as the practice that ``considers harmony not as commanded, but as commanding, and not as a servant, but as the mistress of the oration.'' In this style, the rigorous rules of counterpoint and harmony were paramount.
\item The \textit{seconda pratica} was the new, expressive style pioneered by Monteverdi. Its guiding principle was encapsulated in his own motto: ``Let the word be the mistress of the harmony and not its slave.'' This approach allowed for a much freer treatment of counterpoint and dissonance to better convey the emotional content of the text.
\end{enumerate}
One critic of this new style was the theorist Giovanni Artusi, who condemned Monteverdi's harmonic innovations. Artusi specifically cited the madrigal \textit{Oh, Mirtillo}, which he heard before its publication, as an example of compositional failure. He found the key changes—from $B$-flat to $B$-natural and back—to be an ``aberrant'' and ``puerile error.'' This debate formed the central tension that Monteverdi navigated throughout his career, as he sought greater expressive freedom by systematically expanding the boundaries of accepted musical logic.

\section{Mazzola's Model of Counterpoint}

Traditional counterpoint theory, as codified by Johann Joseph Fux in his 1725 treatise \textit{Gradus ad Parnassum} \cite{fux1965study}, relies on a set of prescriptive rules governing the movement between consonant intervals. This section introduces a model developed by Guerino Mazzola \cite{gM02,AJM15} based on symmetries, which recasts these rules within a formal system based on algebraic symmetries. This model allows for a more precise and quantitative analysis of contrapuntal motion, providing a powerful tool for examining the logic behind Monteverdi's techniques.

The twelve tones of the equal-tempered scale are represented by the ring of integers modulo 12, $\mathbb{Z}/12\mathbb{Z}$. An interval is the distance between two tones, so it is also an element of $\mathbb{Z}/12\mathbb{Z}$. Here it is important to mention that while equal temperament is often considered a later development with respect to the early 17th century, recent scholarship has shown that it was already being advocated and practically implemented in instrumental music during Monteverdi's time \cite{10.1093/em/caw011}. In any case, although the tuning systems of Monteverdi's era varied, it is possible to accommodate the interpretation of his music within an equally (or at least well) tempered framework.

The twelve possible intervals are partitioned into a set of consonances,
\[
 K = \{0, 3, 4, 7, 8, 9\}
\]
which correspond to unison, minor/major thirds, perfect fifth, minor/major sixths, respectively, and the complementary set of dissonances,
\[
 D = \{1, 2, 5, 6, 10, 11\}.
\]

This partition of $\mathbb{Z}/12\mathbb{Z}$ is a \emph{strong dichotomy}, meaning there exists a unique affine transformation, called a \emph{polarity}, that maps the set of consonances to the set of dissonances and vice-versa. For the standard set of consonances $K$, this unique polarity is the function
\begin{align*}
p:\mathbb{Z}/12\mathbb{Z}&\to \mathbb{Z}/12\mathbb{Z}\\
 x&\mapsto 5x + 2.
\end{align*}

This algebraic property establishes a fundamental, symmetric relationship between consonance and dissonance.

The model extends these concepts to analyse the movement between intervals in two-voice counterpoint. Thus two-voice counterpoint is modelled using dual numbers over $\mathbb{Z}/12\mathbb{Z}$, denoted by $\mathbb{Z}/12\mathbb{Z}[\epsilon]$. A dual number $x + \epsilon.y\in\mathbb{Z}/12\mathbb{Z}[\epsilon]$ represents a counterpoint interval where $y$ is the interval of the upper voice (\textit{discantus}) over a fixed voice (\textit{cantus firmus}) pitch $x$.

The set of all consonant intervals is thus
\[
 K[\epsilon] = \{x + \epsilon.k : x \in \mathbb{Z}/12\mathbb{Z}, k \in K\}.
\]

Now we introduce \emph{symmetries} as affine transformations of the space $\mathbb{Z}/12\mathbb{Z}[\epsilon]$, i.e.,
\begin{align*}
g:\mathbb{Z}/12\mathbb{Z}[\epsilon]&\to \mathbb{Z}/12\mathbb{Z}[\epsilon],\\
(x+\epsilon.y)&\mapsto (u_{1}+\epsilon.v_{1})(x+\epsilon.y)+(u_{2}+\epsilon.v_{2})
\end{align*}
where $u_{1}+\epsilon.v_{1}\in \mathbb{Z}/12\mathbb{Z}[\epsilon]^{\times}$ and $u_{2}+\epsilon.v_{2}\in \mathbb{Z}/12\mathbb{Z}[\epsilon]$. In particular, there exists a symmetry which is an \textit{induced polarity} for the dichotomy $(K[\epsilon]/D[\epsilon])$, namely
\[
e^{\epsilon.2}\circ 5.
\]

Thus we can define \textit{counterpoint symmetries} for a consonant interval $\xi\in K[\epsilon]$ as those $g$ that satisfy three conditions:
\begin{enumerate}
\item $\xi\in K[\epsilon]\cap g(D[\epsilon])$,
\item $p$ is a polarity of $(gK[\epsilon]/gD[\epsilon])$ and
\item $g$ maximizes the cardinality of $gK[\epsilon]\cap K[\epsilon]$ among the symmetries that satisfy the first two conditions.
\end{enumerate}

This models the musical idea that in a specific compositional moment, not all consonant moves are equally desirable; the symmetry privileges a certain subset of moves, thereby formalizing the composer's momentary intent or style.

An \textit{admissible successor} is a consonant interval $\eta$ that is a permitted next step from a given consonant interval $\xi$, according to a specific counterpoint symmetry $g$.

It is important to note that the symmetry $g$ can always be chosen among the subgroup $H$ of symmetries of the form
\[
e^{\epsilon.\mathbb{Z}/12\mathbb{Z}}\circ \mathbb{Z}/12\mathbb{Z}[\epsilon]^{\times}.
\]

Thus we use this group as the natural setting to gauge the number of symmetries that mediate in admissible steps. Furthermore, the Hichert algorithm yields the specific symmetries and the number of admitted successors for every consonance. The \textit{Kleiner Kontrapunktsatz} states that any consonance has at least $42$ admitted successors \cite{AJM15}.

\section{Mathematical Models of Tonal Modulation}

Modulation, the process of changing from one key (or \textit{tonality}) to another, is fundamental to creating large-scale musical structure and narrative. It is the engine of harmonic development in Western music. This section presents two complementary mathematical models that formalize the principles of tonal modulation, building upon foundational concepts first explored by theorists like Jean-Philippe Rameau and Arnold Schönberg \cite{schoe1922}.
\subsection{Mazzola's Quantum Model of Modulation}
Developed by Guerino Mazzola, the quantum model formalizes the process of moving between tonalities as a discrete, quantifiable event.

A \textit{tonality} is defined as a covering of the pitch-class space $\mathbb{Z}/12\mathbb{Z}$ by the seven triads (three-note chords) built on the degrees of a musical scale. For example, the C major tonality consists of the seven chords
\[
 \{C, dm, em, F, G, am, B^{\circ}\}.
\]

For a general tonality, we will denote its degrees with Roman numerals, i.e.,
\[
\{I,II,III,IV,V,VI,VII\}.
\]

A \textit{cadential set} is a set of chords that uniquely identifies a tonality. If a set of chords appears in only one tonality and no smaller subset has this property, it is a \textit{cadence}.

There are five such minimal sets for the major scale, namely,
\[
k_{1} = \{II, V\}, k_{2} = \{II,III\}, k_{3}=\{III,IV\}, k_{4}=\{IV,V\}, k_{5}=\{VII\}.
\]

A \textit{modulator} is a symmetry $m$ that maps the starting tonality $S$ to the target tonality $T$ (i.e., $m(S) = T$).

A \textit{modulation quantum} $M$ is the minimal set of notes required to execute the modulation. It is the minimal set that satisfies three conditions:
\begin{enumerate}
\item It is invariant under the modulator $m$ (i.e. $m(M) = M$).
\item It contains all the notes of the cadence $k$ in the target tonality.
\item Its intersection $m\cap T$ with the target tonality $T$ is rigid (it has no non-trivial symmetries). 
\end{enumerate}

The chords common to both the starting and target tonalities that fall within this quantum $M$ are the \textit{pivot chords} that enable a smooth transition.

\begin{example}\label{ejemplo3_3}
For the modulation from $C$ to $D$ with cadence $k_{4}=\{IV_{D},V_{D}\}$, we have the modulator $m=e^{6}\cdot-1$. The modulation quantum is
\[
M=\{1,2,4,5,7,9,11\}=II_{C}\cup III_{C}\cup V_{C}\cup VII_{C}\cup \underbrace{II_{D}\cup VII_{D}\cup IV_{D} \cup V_{D}}_{\text{Pivots}}.
\]
Indeed, since $C=\{0,2,4,5,7,9,11\}$ and $D=\{2,4,6,7,9,11,1\}$ we have $M\cap D=\{1,2,4,7,9,11\}$, which we know is rigid because it is precisely a strong dichotomy.
\end{example}

\begin{example}\label{ejemplo3_4}
The major scale with tonic $G$ is
\[
G=\{7,9,11,0,2,4,6\}.
\]

For the modulation from $C$ to $G$ with cadence $k_{5}=\{VII_{G}\}$, we have the modulator $m=e^{11}\cdot-1$. The modulation quantum is
\[
M=\{0,2,5,6,9,11\}=II_{C}\cup IV_{C}\cup VII_{C}\cup \underbrace{III_{G}\cup V_{G}\cup VII_{G}}_{\text{Pivots}}.
\]
We have $M\cap G=\{0,2,6,9,11\}$ and it is easy to check it is rigid. The scale degrees from $M$ that are also in $G$ are the pivots:
\[
III_{G}=\{11,2,6\}, \quad V_{G}=\{2,6,9\}, \quad VII_{G}=\{6,9,0\}.
\]
Since $e^{11}\cdot-1(\{0,2,6,9,11\})=\{11,9,5,2,0\}$, we must add $5$ to $M$ to ensure that $m$ is its symmetry, because $e^{11}\cdot-1(5)=-5+11=6$. This proves that $M$ is minimal.
\end{example}

\begin{example}\label{ejemplo3_5}
The major scale with tonic $E\flat$ is
\[
E\flat=\{3,5,7,8,10,0,2\}.
\]

For the modulation from $C$ to $E\flat$ with cadence $k_{1}=\{II_{E\flat},V_{E\flat}\}$, we have the modulator $m=e^{7}\cdot11$. The modulation quantum is
\[
M = \{0,2,5,7,8,9,10,11\} = II_C \cup V_C \cup VII_C \cup \underbrace{II_{E\flat} \cup III_{E\flat} \cup V_{E\flat} \cup VII_{E\flat}}_{\text{Pivots}}
\]
The intersection is
\[
M\cap E\flat = \{0,2,5,7,8,10\},
\]
and it is rigid because it is a strong dichotomy. The scale degrees from $M$ that are also in $E\flat$ are the pivots:
\[
II_{E\flat}, \quad III_{E\flat}, \quad V_{E\flat}, \quad VII_{E\flat}.
\]

Note that we can take $E\flat$ as equivalent to $D\sharp$ because of enharmonic identification.
\end{example}

Let us summarise the information from the three previous examples in a table.

\[
\begin{tabular}{|c|c|c|c|c|}
\hline
\textbf{Modulation} & $m$ & \textbf{Cadence} & \textbf{$M\cap T$} & \textbf{Pivots} \\
\hline
$C\to D$ & $e^6\cdot-1$ & $k_4=\{IV_D,V_D\}$ & $\{1,2,4,7,9,11\}$ & $II_D, IV_D, V_D, VII_D$ \\
\hline
$C\to G$ & $e^{11}\cdot-1$ & $k_5=\{VII_G\}$ & $\{0,2,6,9,11\}$ & $III_G, V_G, VII_G$ \\
\hline
$C\to E\flat$ & $e^7\cdot11$ & $k_1=\{II_{E\flat},V_{E\flat}\}$ & $\{0,2,5,7,8,10\}$ & $II_{E\flat}, III_{E\flat}, V_{E\flat}, VII_{E\flat}$ \\
\hline
\end{tabular}
\]

\subsection{The Duality Model of Modulation}

A problem with Mazzola's model for modulation is that it does not account for transitions between minor and major tonalities. Here 
the group-theoretic $PLR/TI$ model describes the elegant structural relationships between major and minor keys.
The model\footnote{For an exposition of the model, see \cite{AAdPLPM09}.} operates on the set of all $24$ major and minor triads, denoted by $M$.

The TI group is the group of affine transformations on $\mathbb{Z}/12\mathbb{Z}$, consisting of \textit{transpositions}
\[
 T_{n}(x) = x + n,
\]
which shift a chord up or down, and \textit{inversions}
\[
 I_{n}(x) = -x + n,
\]
which flip a chord around an axis.

The PLR group, on the other hand, consists of three fundamental musical transformations that relate major and minor chords:
\begin{enumerate}
\item $P$ (Parallel): Transforms a major chord to its parallel minor (e.g., $F$ major to $f$ minor) or vice-versa.
\item $L$ (\textit{Leittonwechsel}): Transforms a major chord to the minor chord on its leading tone (e.g., $C$ major to $e$ minor).
\item $R$ (Relative): Transforms a major chord to its relative minor (e.g., $C$ major to $a$ minor) or vice-versa.
\end{enumerate}

A crucial result is that the algebraic group TI and the music-theoretic group PLR are isomorphic. Furthermore, their actions on the set of $24$ triads $M$ are dual and regular. This mathematical property proves a deep, non-obvious structural equivalence between geometric pitch-space operations (transposition, inversion) and functional harmonic relationships (parallel, relative), a connection composers like Monteverdi intuitively grasped. This means that any chord can be transformed into any other chord (transitivity), the transformation is unique (freeness), and the elements of the two groups commute.

In particular, we can consider the cadence $k_{4}=\{IV,V\}$ as defining the major tonality, with relations $T^{7}(I)=V$ and $T^{5}(I)=IV$, and thus call it the \emph{major} cadence. Using the fact that $T$ commutes with $R$, we can recover the $k_{2}$ cadence since
\begin{align}
R\circ T^{7}(I)&=T^{7}\circ R(I)=T^{7}(VI)=III,\\
R\circ T^{5}(I)&=T^{5}\circ R(I)=T^{5}(VI)=II.
\end{align}

Hence $k_{2}$ can be regarded as a \emph{minor} cadence, and the modulation is accomplished by displaying it after the major cadence\footnote{Rimsky-Korsakoff, for example, describes a modulatory process very similar to this from $C$ major to $a$ minor \cite[p.~44]{rimsky1910harmonie}.}. Of course, the converse modulation proceeds analogously.

\section{Application and Analysis of Monteverdi's Works}
This section forms the analytical core of the paper, applying the previously outlined mathematical models to musical excerpts from Monteverdi's sacred works and his landmark opera \textit{L'Orfeo}. These analyses will demonstrate the explanatory power of the models in formalizing the structural logic behind Monteverdi's compositional choices.

\subsection{Contrapuntal Analysis via Symmetries}
We first apply the symmetry-based counterpoint model to several excerpts featuring common harmonic patterns in Monteverdi's music.

\subsubsection{Analysis of Progressions of Fifths and Thirds}

We will examine certain common transitions in Monteverdi’s works, valid both in the counterpoint model and in Fux’s rigorous one. The fragments we analyse are some of those that appear in an article by Kang \cite{kang2011monteverdi}. This first group corresponds to what Kang calls ``$5-3$ harmonies'', that is, alternating intervals of a fifth and a third.

\begin{figure}[ht]
    \centering
    \includegraphics[width=\linewidth]{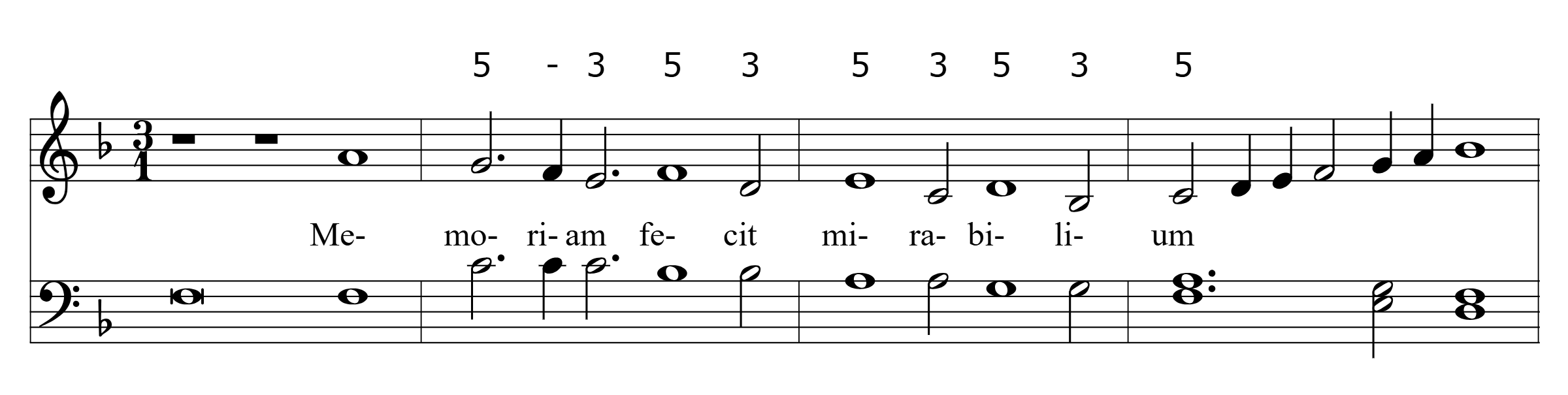}
    \caption{\textit{Confitebor primo (1640)}, measures from 64 to 67}
    \label{confitebor}
\end{figure}

The first is \textit{Confitebor primo (1640)} \cite{cM40}, in measures 64–67. The following succession of counterpoint intervals is obtained

\begin{multline*}
\xi_{1}=0+\epsilon.7,\xi_{2}=0+\epsilon.4,\xi_{3}=10+\epsilon.7,\xi_{4}=10+\epsilon.4,\\
\xi_{5}=9+\epsilon.7,\xi_{6}=9+\epsilon.3,\xi_{7}=7+\epsilon.7,\xi_{8}=7+\epsilon.4,\xi_{9}=5+\epsilon.7.
\end{multline*}

The counterpoint symmetry sets associated with this sequence are
\[
\begin{tabular}{c c }
Set of symmetries & Cardinality of the set \\
$\{g_{1}=e^{\epsilon.0}\circ7\}$ & 1 \\
$\{g_{2}=e^{\epsilon.6}\circ(7+\epsilon.6)\}$ & 1 \\
$\{g_{3}=e^{\epsilon.0}\circ7\}$ & 1 \\
$\{g_{4}=e^{\epsilon.6}\circ(1+\epsilon.6)\}$ & 1 \\
$\{g_{5}=e^{\epsilon.0}\circ7\}$ & 1 \\
$\{g_{6}=e^{\epsilon.8}\circ(5+\epsilon.4)\}$ & 1 \\
$\{g_{7}=e^{\epsilon.0}\circ7\}$ & 1 \\
$\{g_{8}=e^{\epsilon.6}\circ(7+\epsilon.6)\}$ & 1
\end{tabular}
\]

\begin{figure}[ht]
    \centering
    \includegraphics[width=\linewidth]{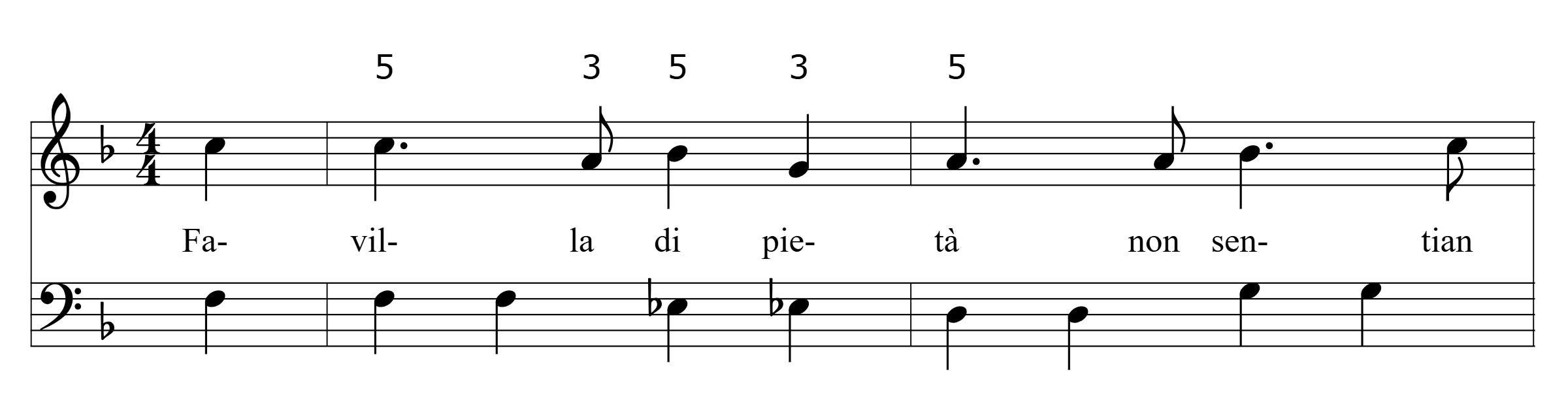}
    \caption{\textit{Ma tu, più che mai dura (1640)}, measures from 6 to 8}
    \label{matu}
\end{figure}

Next we have a fragment from \textit{Ma tu, più che mai dura}, from Book V of madrigals \cite{cM08}, in measures 6–8, yielding the following counterpoint intervals
\[
\xi_{1}=5+\epsilon.7,\xi_{2}=5+\epsilon.4,\xi_{3}=3+\epsilon.7,\xi_{4}=3+\epsilon.4,\xi_{5}=2+\epsilon.7.
\]
The counterpoint symmetry sets are
\[
\begin{tabular}{c c }
Set of symmetries & Cardinality of the set \\
$\{g_{1}=e^{\epsilon.0}\circ7\}$ & 1 \\
$\{g_{2}=e^{\epsilon.6}\circ(7+\epsilon.6)\}$ & 1 \\
$\{g_{3}=e^{\epsilon.0}\circ7\}$ & 1 \\
$\{g_{4}=e^{\epsilon.6}\circ(1+\epsilon.6)\}$ & 1
\end{tabular}
\]

\begin{figure}[ht]
    \centering
    \includegraphics[width=\linewidth]{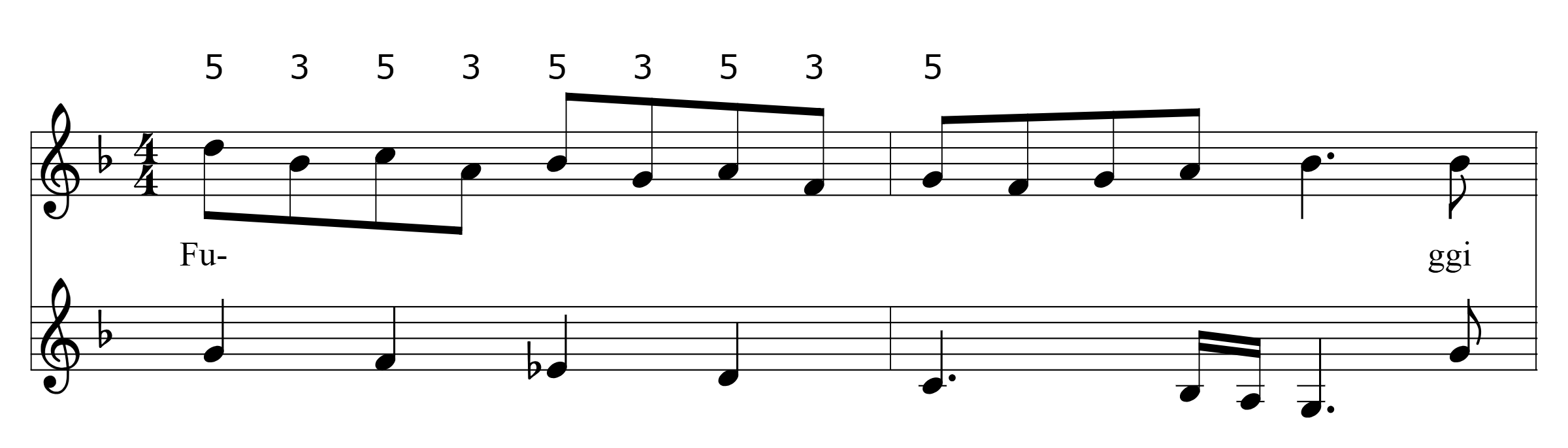}
    \caption{\textit{Io mi son giovinetta}, measures from 52 to 53}
    \label{iomi}
\end{figure}

Finally we have a fragment from \emph{Io mi son giovinetta}, from Book IV of madrigals \cite{cM03}, in measures 52–53, yielding the following counterpoint intervals
\begin{multline*}
\xi_{1}=7+\epsilon.7,\xi_{2}=7+\epsilon.3,\xi_{3}=5+\epsilon.7,\xi_{4}=5+\epsilon.4,\\
\xi_{5}=3+\epsilon.7,\xi_{6}=3+\epsilon.4,\xi_{7}=2+\epsilon.7,\xi_{8}=2+\epsilon.3,\xi_{9}=0+\epsilon.7.
\end{multline*}
The counterpoint symmetry sets are
\[
\begin{tabular}{c c }
Set of symmetries & Cardinality of the set \\
$\{g_{1}=e^{\epsilon.0}\circ7\}$ & 1 \\
$\{g_{2}=e^{\epsilon.8}\circ(5+\epsilon.4)\}$ & 1 \\
$\{g_{3}=e^{\epsilon.0}\circ7\}$ & 1 \\
$\{g_{4}=e^{\epsilon.6}\circ(7+\epsilon.6)\}$ & 1 \\
$\{g_{5}=e^{\epsilon.0}\circ7\}$ & 1 \\
$\{g_{6}=e^{\epsilon.6}\circ(1+\epsilon.6)\}$ & 1 \\
$\{g_{7}=e^{\epsilon.0}\circ7\}$ & 1 \\
$\{g_{8}=e^{\epsilon.8}\circ(5+\epsilon.4)\}$ & 1
\end{tabular}
\]

\subsubsection{Analysis of Octaves and Sixths}

This corresponds to what Kang calls $8$--$6$ harmonies, i.e., alternations of octaves and sixths. The first fragment of this type comes from the \emph{Laudate Dominum (1640)} \cite{cM40}, measures 98--103, with the following intervals.

\begin{figure}[ht]
    \centering
    \includegraphics[width=\linewidth]{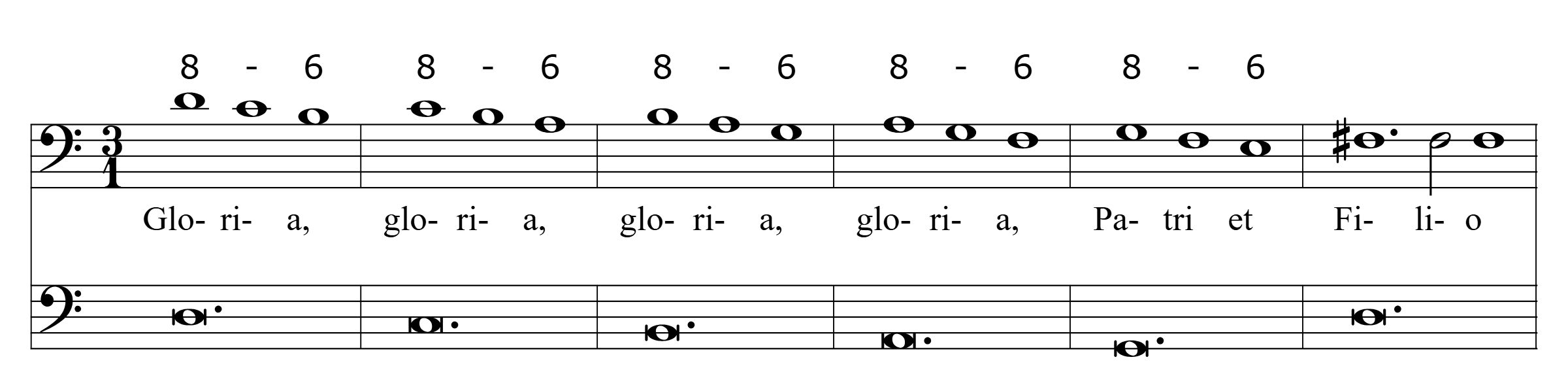}
    \caption{\textit{Laudate Dominum (1640)}, measures from 98 to 103}
    \label{laudate}
\end{figure}

\begin{multline*}
    \xi_{1}=2+\epsilon.0,\xi_{2}=2+\epsilon.9,\xi_{3}=0+\epsilon.0,\xi_{4}=0+\epsilon.9,\xi_{5}=11+\epsilon.0,\\
    \xi_{6}=11+\epsilon.8,\xi_{7}=9+\epsilon.0,\xi_{8}=9+\epsilon.8,\xi_{9}=7+\epsilon.0,\xi_{10}=7+\epsilon.9\\
\end{multline*}

The counterpoint symmetries are as follows.

\begin{small}
\[
\begin{tabular}{c c }
    Set of symmetries & Cardinality \\[-0.5em]
     & of the set\\
    $\{g_{1}=e^{\epsilon.6}\circ(1+\epsilon.6), g_{2}=e^{\epsilon.6}\circ(7+\epsilon.6)\}$ & 2 \\
    $\{g_{3}=e^{\epsilon.8}\circ(5+\epsilon.8), g_{4}=e^{\epsilon.8}\circ(5+\epsilon.4)\}$ & 2 \\
    $\{g_{5}=e^{\epsilon.6}\circ(1+\epsilon.6), g_{6}=e^{\epsilon.6}\circ(7+\epsilon.6)\}$ & 2 \\
    $\{g_{7}=e^{\epsilon.8}\circ(5+\epsilon.8), g_{8}=e^{\epsilon.8}\circ(5+\epsilon.4)\}$ & 2 \\
    $\{g_{9}=e^{\epsilon.11}\circ(11+\epsilon.8), g_{10}=e^{\epsilon.11}\circ(11+\epsilon.4), g_{11}=e^{\epsilon.11}\circ11\}$ & 3 \\
    $\{g_{12}=e^{\epsilon.3}\circ(7+\epsilon.8), g_{13}=e^{\epsilon.3}\circ7\}$ & 2 \\
    $\{g_{14}=e^{\epsilon.11}\circ(11+\epsilon.8), g_{15}=e^{\epsilon.11}\circ(11+\epsilon.4), g_{16}=e^{\epsilon.11}\circ11\}$ & 3 \\
    $\{g_{17}=e^{\epsilon.3}\circ(7+\epsilon.8), g_{18}=e^{\epsilon.3}\circ7\}$ & 2 \\
    $\{g_{19}=e^{\epsilon.6}\circ(1+\epsilon.6), g_{20}=e^{\epsilon.6}\circ(7+\epsilon.6)\}$ & 2 \\
\end{tabular}
\]
\end{small}

\begin{figure}[ht]
    \centering
    \includegraphics[width=\linewidth]{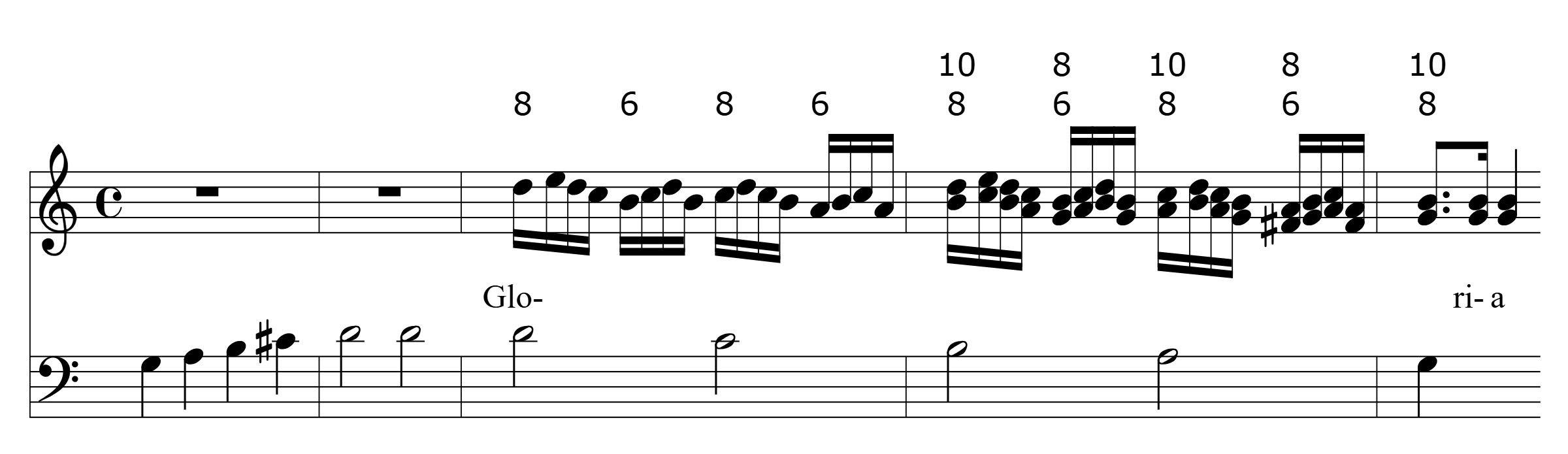}
    \caption{\textit{Gloria a 7}, measures from 1 to 5}
    \label{gloriaa7}
\end{figure}

We follow with an excerpt from \emph{Gloria a 7 (1640)} \cite{cM40}, in measures 1–5. This fragment is interesting because the counterpoint appears starting from the third measure, and the successions $8$–$6$ and $10$–$8$ arise from considering different pairs of voices. The following counterpoint intervals are obtained for the $8$–$6$ succession:

\begin{multline*}
\xi_{1}=2+\epsilon.0,\xi_{2}=2+\epsilon.9,\xi_{3}=0+\epsilon.0,\xi_{4}=0+\epsilon.9,\\
\xi_{5}=11+\epsilon.0,\xi_{6}=11+\epsilon.8,\xi_{7}=9+\epsilon.0,\xi_{8}=9+\epsilon.8,\xi_{9}=7+\epsilon.0.
\end{multline*}

The counterpoint symmetries associated with the succession are:
\[
\begin{tabular}{c c }
Set of symmetries & Cardinality\\
  & of the set\\\
$\{g_{1}=e^{\epsilon.6}\circ(1+\epsilon.6), g_{2}=e^{\epsilon.6}\circ(7+\epsilon.6)\}$ & 2 \\
$\{g_{3}=e^{\epsilon.8}\circ(5+\epsilon.8), g_{4}=e^{\epsilon.8}\circ(5+\epsilon.4)\}$ & 2 \\
$\{g_{5}=e^{\epsilon.6}\circ(1+\epsilon.6), g_{6}=e^{\epsilon.6}\circ(7+\epsilon.6)\}$ & 2 \\
$\{g_{7}=e^{\epsilon.8}\circ(5+\epsilon.8), g_{8}=e^{\epsilon.8}\circ(5+\epsilon.4)\}$ & 2 \\
$\{g_{9}=e^{\epsilon.11}\circ(11+\epsilon.8), g_{10}=e^{\epsilon.11}\circ(11+\epsilon.4), g_{11}=e^{\epsilon.11}\circ11\}$ & 3 \\
$\{g_{12}=e^{\epsilon.3}\circ7, g_{13}=e^{\epsilon.3}\circ(7+\epsilon.8)\}$ & 2 \\
$\{g_{14}=e^{\epsilon.11}\circ(11+\epsilon.8), g_{15}=e^{\epsilon.11}\circ(11+\epsilon.4), g_{16}=e^{\epsilon.11}\circ11\}$ & 3 \\
$\{g_{17}=e^{\epsilon.3}\circ7, g_{18}=e^{\epsilon.3}\circ(7+\epsilon.8)\}$ & 2
\end{tabular}
\]

The counterpoint intervals for the $10$–$8$ succession are:
\[
\xi_{1}=11+\epsilon.3,\xi_{2}=11+\epsilon.0,\xi_{3}=9+\epsilon.3,\xi_{4}=9+\epsilon.0,\xi_{5}=7+\epsilon.4.
\]

The counterpoint symmetries associated with this succession are:
\[
\begin{tabular}{c c }
Set of & Cardinality of the \\
symmetries & set\\
$\{g_{1}=e^{\epsilon.8}\circ(5+\epsilon.8), g_{2}=e^{\epsilon.8}\circ(5+\epsilon.4)\}$ & 2 \\
$\{g_{3}=e^{\epsilon.6}\circ(1+\epsilon.6), g_{4}=e^{\epsilon.6}\circ(7+\epsilon.6),$ & \multirow{2}{0.5em}{5}\\
$g_{5}=e^{\epsilon.11}\circ(11+\epsilon.8),
g_{6}=e^{\epsilon.11}\circ(11+\epsilon.4),
g_{7}=e^{\epsilon.11}\circ11\}$\\
$\{g_{8}=e^{\epsilon.8}\circ(5+\epsilon.8), g_{9}=e^{\epsilon.8}\circ(5+\epsilon.4)\}$ & 2 \\
$\{g_{10}=e^{\epsilon.11}\circ(11+\epsilon.4), g_{11}=e^{\epsilon.11}\circ11\}$ & 2
\end{tabular}
\]

Kang thus identifies several passages across Monteverdi's sacred and secular output where successions of 5-3, 8-6 or 8-10 intervallic patterns serve both structural and expressive functions. In her analysis of \textit{Confitebor primo}, she notes how Monteverdi avoids consecutive fifths between alto and bass by having the upper voices leap between pitches of a descending line, creating successive 5-3 intervals. Similarly, she observes that in \textit{Ma tu, più che mai dura}, the juxtaposition of 5-3 harmonies on $E\flat$ and $D$ emphasizes the word \textit{pietà}, while in \textit{Io mi son giovinetta}, rapid 5-3 successions create a madrigalistic effect on \textit{fuggi}. Kang treats these as instances of clever contrapuntal technique—sometimes serving textual expression, sometimes functioning as independent compositional devices.

Our group-theoretic analysis reveals a more systematic pattern underlying these observations. The 5-3 successions Kang identifies are indeed valid under both Mazzola's counterpoint model and traditional Fuxian rules (though they would constitute ``hidden'' parallels in second-species counterpoint). What makes them structurally distinctive is their parsimony: each transition is mediated by a single counterpoint symmetry. 

More significantly, in all three 5-3 examples, the transformation $g_4 = e^{\epsilon.6}\circ(1+\epsilon.6)$ (the only symmetry along the passages that fixes the \textit{cantus firmus}) occurs precisely at the expressively charged words. In \textit{Confitebor primo}, this transformation marks \textit{mirabilium} (wonders), a theologically significant term that Kang does not specifically highlight. In \textit{Ma tu, più che mai dura}, it emphasizes again \textit{pietà} (compassion), and in \textit{Io mi son giovinetta}, it coincides with the chromatic $E\flat$ on \textit{fuggi} (flee), literally a note that ``flees'' from the key signature.

This pattern suggests that what Kang observes as varied applications of a contrapuntal technique (``sometimes for textual expression, sometimes as an independent device'') may actually reflect a consistent compositional strategy: Monteverdi systematically deploys the \textit{cantus-firmus}-fixing transformation to mark words of emotional or theological weight.

Kang also examines 8-6 progressions in \textit{Laudate Dominum} and \textit{Gloria a 7}, where consecutive octaves are avoided by descending to the sixth. Our analysis shows these differ structurally from the 5-3 examples: each transition requires at least two symmetries, making them less parsimonious. In \textit{Laudate Dominum}, transitions from perfect to imperfect consonances involve three simultaneous symmetries: $\{e^{\epsilon.11}\circ(11+\epsilon.8), e^{\epsilon.11}\circ(11+\epsilon.4), e^{\epsilon.11}\circ11\}$.

The \textit{Gloria a 7} is particularly revealing. Beyond the 8-6 progressions, it contains 10-8 harmonies (which can be thought as 8-3: octave to third) in the fourth measure. One transition here exhibits five simultaneous symmetries: $\{e^{\epsilon.6}\circ(1+\epsilon.6), e^{\epsilon.6}\circ(7+\epsilon.6), e^{\epsilon.11}\circ(11+\epsilon.8), e^{\epsilon.11}\circ(11+\epsilon.4), e^{\epsilon.11}\circ11\}$. These jumps in cardinality seem to correlate with the entrance of additional voices in the homophonic \textit{tutti} sections, precisely the moments Kang identifies as ``especially interesting.'' The mathematical structure thus reflects the textural complexity: more voices require more simultaneous symmetries, creating less parsimonious but texturally richer passages.

In summary, our analysis both validates and extends Kang's observations. The 5-3 and 8-6 patterns she identifies are indeed compositionally significant, but their significance lies not merely in avoiding parallels or creating local effects—rather, they reflect systematic structural principles that Monteverdi deployed consistently across nearly four decades of compositional practice.

\subsection{Modulation Analysis in L'Orfeo}
We will focus the modulations present in different fragments of the opera \textit{L'Orfeo}, which validate the quantum model of modulation and modulation between major and minor. The fragments of our analysis were taken from the article \textit{From Modal to Tonal: The Influence of Monteverdi on Musical Development} \cite{Perritt2017FromMT}.

\begin{figure}[p]
\centering \begin{overpic}[width=\linewidth]{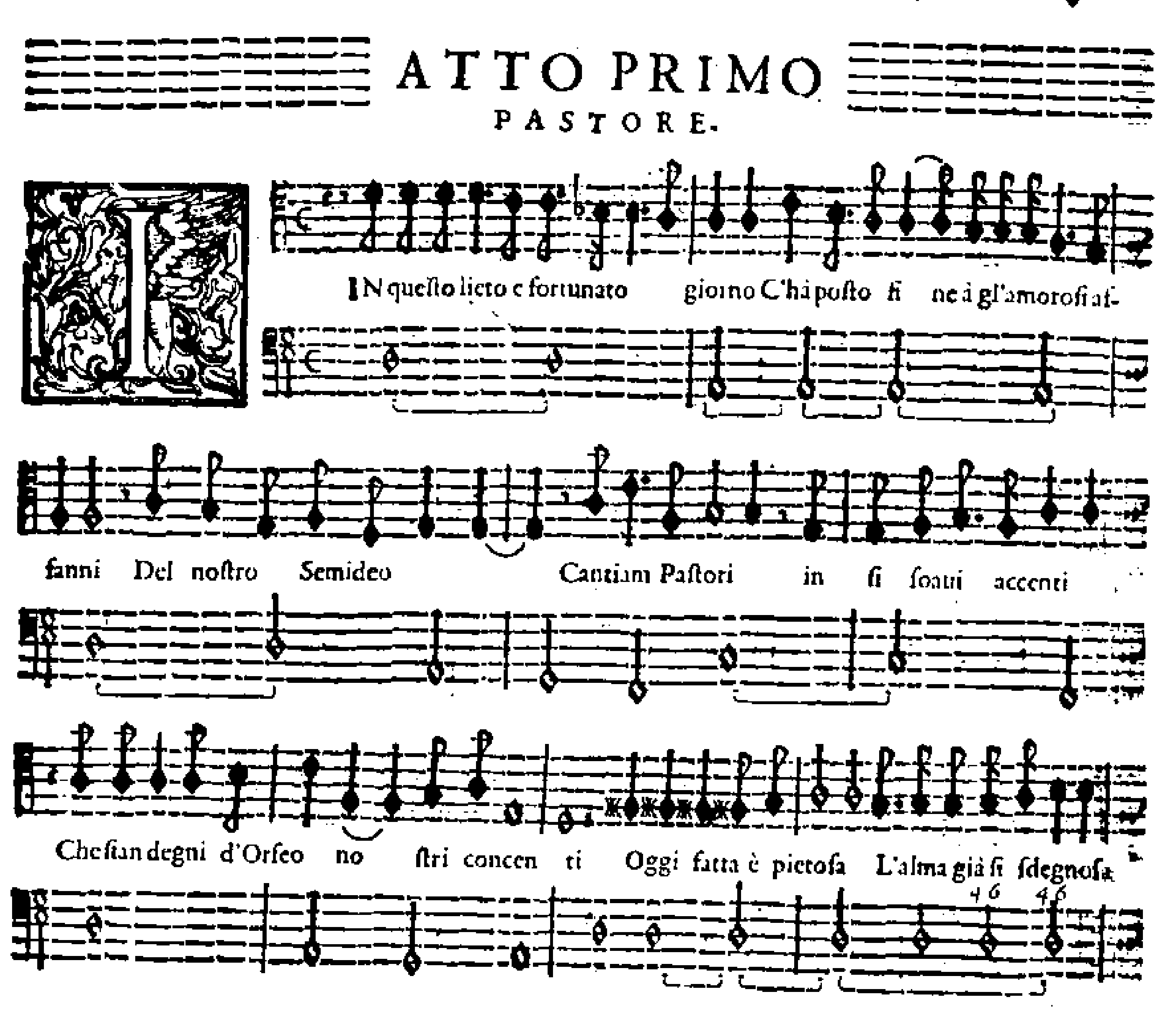}
 \put (38,25.5) {\insquare{$k_{4}=\{F,\hspace{1.7em} G\}$}}
 \put (27,1) {\insquare{$k_{4}=\{G,\qquad A\}$}}
\end{overpic}    
    \caption{\textit{In questo lieto e fortunato giorno}, first aria of Act I of \textit{L'Orfeo}, measures 1 to 8}
    \label{in_questo1}

\end{figure}

\begin{figure}[p]
\centering \begin{overpic}[width=\linewidth]{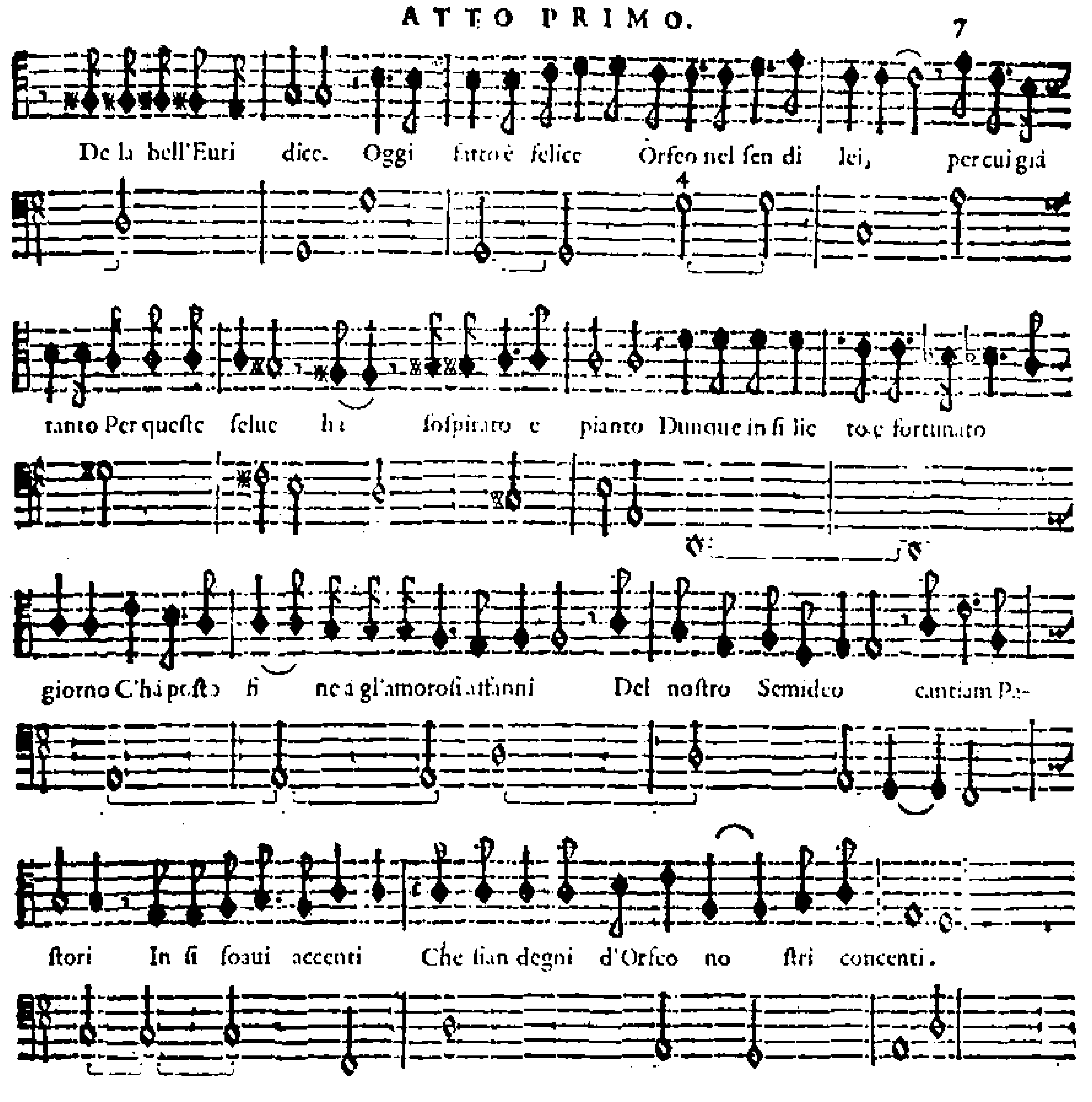}
 \put (0,49) {\insquare{$k_{5}=\{g\sharp^{\circ}\}$}}
 \put (55,24) {\insquare{$k_{1}=\{D,\hspace{5em} G\}$}}
 \end{overpic}
    \caption{\textit{In questo lieto e fortunato giorno}, first aria of Act I of \textit{L'Orfeo}, measures 9 to 20}
    \label{in_questo2}
\end{figure}

We begin by examining the modulations present in the aria \textit{In questo lieto e fortunato giorno}, where the intersection of the modulation quantum with the scale of the destination tonality is also a strong dichotomy. This can be seen reproduced in figures \ref{in_questo1} and \ref{in_questo2}.


In this aria, we can find the modulation from $C$ to $D$, which was calculated in Example \ref{ejemplo3_3}. Recall that the chords of $D$ are:
\[
\begin{tabular}{c c}
    $I_{D}=\{2,6,9\}$, &$D$, \\
    $II_{D}=\{4,7,11\}$, &$e$, \\
    $III_{D}=\{6,9,1\}$, &$f\sharp$, \\
    $IV_{D}=\{7,11,2\}$, &$G$, \\
    $V_{D}=\{9,1,4\}$, &$A$, \\
    $VI_{D}=\{11,2,6\}$, &$b$, \\
    $VII_{D}=\{1,4,7\}$, &$c\sharp^{\circ}$.
\end{tabular}
\]

The cadence $k_{4}$ for $C$ is $\{F,G\}$, and it appears in measure $4$ in figure \ref{in_questo1}. In measure $5$ the first degree of $D$ appears. Recall that a cadence for $D$ is $k_{4}=\{G,A\}$ and the pivots are $e$, $G$, $A$ and $c\sharp^{\circ}$ for this modulation. In this case, note that the cadence is a subset of the pivots, and its chords appear in measure $6$.

Now, Monteverdi modulates from $D$ to $A$, which we already calculated in Example \ref{ejemplo3_4}, although it must be transposed up one tone. Recall that the chords of $A$ are:
\[
\begin{tabular}{c c}
    $I_{A}=\{9,1,4\}$, &$A$, \\
    $II_{A}=\{11,2,6\}$, &$b$, \\
    $III_{A}=\{1,4,8\}$, &$c\sharp$, \\
    $IV_{A}=\{2,6,9\}$, &$D$, \\
    $V_{A}=\{4,8,11\}$, &$E$, \\
    $VI_{A}=\{6,9,1\}$, &$f\sharp$, \\
    $VII_{A}=\{8,11,2\}$, &$g\sharp^{\circ}$.
\end{tabular}
\]

The pivots are $c\sharp$, $E$ and $g\sharp^{\circ}$. The cadence $k_{5}$ for $A$ is $\{g\sharp^{\circ}\}$, and the only chord that comprises it appears in measure $12$ (see figure \ref{in_questo2}), taking advantage of the fact that it is also a pivot. The first degree of $A$ appears in the same measure as the cadence.

We continue with the final modulation from $A$ to $C$, which was calculated in Example \ref{ejemplo3_5}, although it must be transposed up three semitones. The cadence $k_{1}$ for $C$ is $\{d,G\}$, whose chords are found in measures $16$ and $17$ and which are also pivots. The first degree of $C$ appears in measure $18$.

To conclude, we will analyse two \textit{ariettas} found in Act II of \textit{L'Orfeo}; the very first is \textit{Ecco pur}. Certain important chords of the \textit{arietta} can be represented with the following commutative diagram:

\begin{figure}[p]
\centering \begin{overpic}[width=\linewidth]{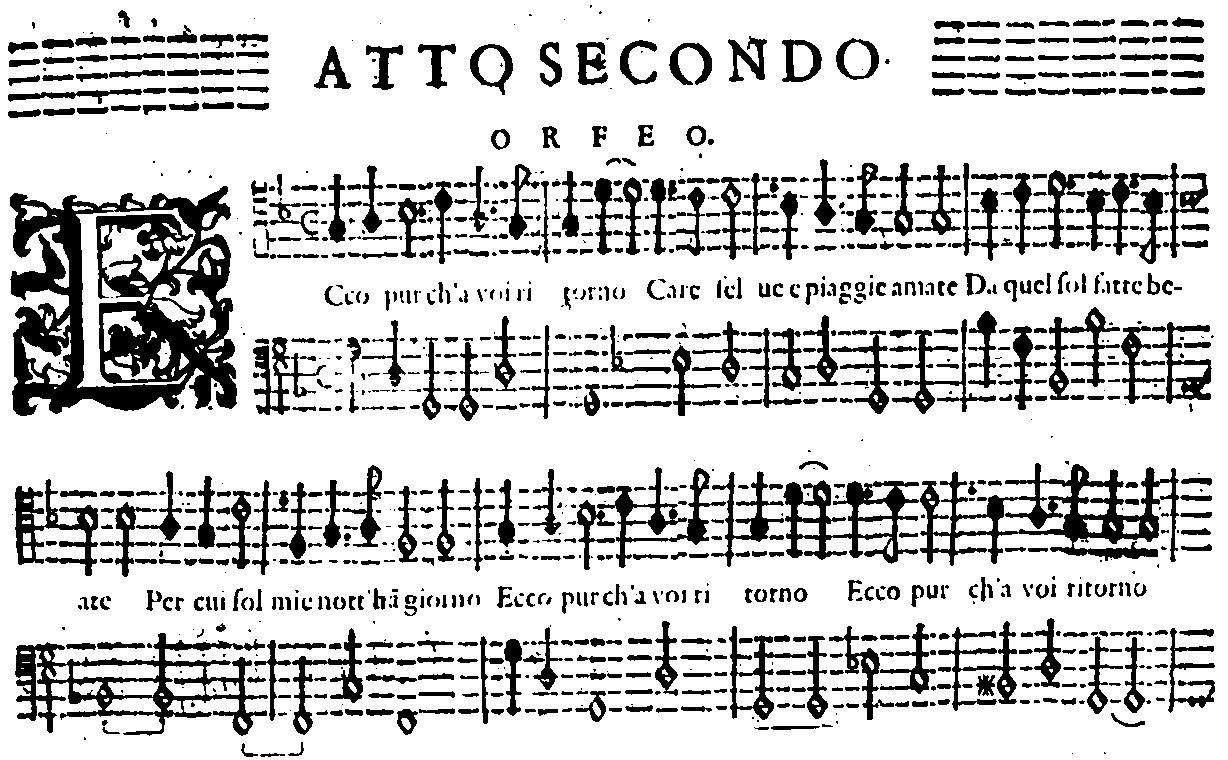}
\put(85,26){\insquare{$B\flat$}}
\put(22,-1.25){\insquare{$F$}}
\put(69,-1.25){\insquare{$E\flat$}}
\put(78,-1.25){\circled{$c$}}
\put(84,-1.25){\circled{$d$}}
\put(89,-1.25){\circled{$g$}}
\end{overpic}
    \caption{Ecco pur, from the second act of L'Orfeo}
    \label{ecco_pur}
\end{figure}

\[ \begin{tikzcd}
\arrow[swap]{d}{R}E\flat & \arrow[swap]{l}{T_5}
B\flat\arrow{r}{T_7} \arrow{d}{R}& F \arrow{d}{R} \\%
c &\arrow{l}{T_5}  g \arrow[swap]{r}{T_7}& d
\end{tikzcd}
\]



Thus for the tonality of $B\flat$, its cadential set is $k_{4}=\{E\flat,F\}$ and for the tonality of $g$ its cadential set is $k_{2}=\{c,d\}$.

From figure \eqref{ecco_pur}, the chords enclosed in squares are major and those enclosed in circles are minor chords. Note that half of the fragment of the aria begins in $B\flat$, where its cadence is slowly displayed, then ends in $c$ together with its cadence very quickly.

The second \textit{arietta} that we will analyse is found in Act II, titled \textit{Mira, deh mira Orfeo}. Here we have the following commutative diagram:

\[ \begin{tikzcd}
\arrow[swap]{d}{R}F & \arrow[swap]{l}{T_5}
C\arrow{r}{T_7} \arrow{d}{R}& G \arrow{d}{R} \\%
d &\arrow{l}{T_5}  a \arrow[swap]{r}{T_7}& e
\end{tikzcd}
\]

\begin{figure}[p]
\centering \begin{overpic}[width=\linewidth]{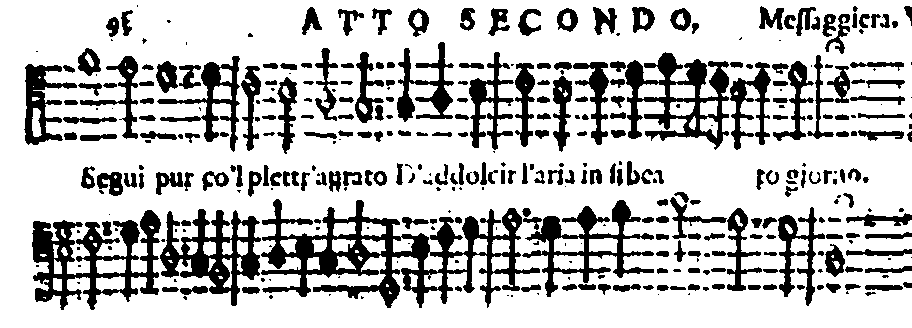}
\put(27,-1.25){\circled{$C$}}
\put(40,-1.25){\circled{$G$}}
\put(7,-1.25){\circled{$F$}}
\put(37.5,-1.25){\insquare{$d$}}
\put(45,-1.25){\insquare{$e$}}
\put(55,-1.25){\insquare{$a$}}
\end{overpic}
    \caption{Mira, deh mira Orfeo, from the second act of L'Orfeo, measures 9 to 12.}
    \label{mira_de_orfeo}
\end{figure}

Hence the cadence $k_{4}=\{F,G\}$ transforms into the cadence $k_{2}=\{d,e\}$. The cadence $k_{2}$ of $a$ minor and the cadence $k_{4}$ of $C$ major appear in measures 10 and 11, respectively, of figure \ref{mira_de_orfeo}.

Perritt examines how Monteverdi navigates the transition from modality to tonality in \textit{L'Orfeo}, highlighting passages where modal and tonal elements coexist. Regarding the aria \textit{In questo lieto e fortunato giorno} from Act I, she notes that the piece begins in Dorian mode on $D$ (with $B\flat$ as the lowered sixth and $C\natural$ as the lowered seventh\footnote{Though it is in fact $G$ Dorian, but this does not properly alters the subsequent analysis.}) but shifts to $D$ major through the introduction of $F\sharp$ and $C\sharp$, supported by strong harmonic progressions. She observes that while initially presented in $C$ major (Ionian mode), the prolongation of the $D$ major chord and $D$ in the bass helps establish a tonal rather than modal framework, calling this ``an interesting combination of modal theory and tonality within a short section.''

Our quantum modulation model provides a more detailed accounting of this aria's tonal structure. Perritt focuses solely on the $C$ to $D$ modulation, but the complete aria presents three modulations, $C \to D \to A \to C$, each corresponding to a distinct dramatic moment:

\begin{itemize}
    \item The call to the shepherds, in $C$.
    \item The beautiful Euridice yields to Orpheus's advances, in $D$. This modulation is not arbitrary. Scholars have noted that Monteverdi systematically associates $D$ major with Euridice throughout \textit{L'Orfeo}. In Orpheus's monologue \textit{Rosa del ciel}, his $G$-Dorian tonal level consistently moves toward $D$ major when Euridice is mentioned: when he first saw her (``pria ti vidi''), when he sighed for her (``per te sospirai''), and finally at the end when he turns ``in an almost scenic gesture, to Euridice's final of D major in anticipation of her response'' \cite{Steinheuer_2007}. Thus the modulation to $D$ in our aria participates in a larger tonal-dramatic strategy spanning the entire opera.
    \item Orpheus's suffering in the woods for love of Euridice, in $A$. This resonates with Monteverdi's association of A major with ``the sphere of the underworld'' \cite{Steinheuer_2007}, foreshadowing his coming journey: suffering for lost love already invokes the underworld's domain.
\end{itemize}

The succession of keys thus narrates Orpheus's emotional journey, culminating in his return to $C$ major when Euridice reciprocates his love—the joy represented by the opening tonality restored (and, of course, Orpheus's return to his departing point at Thrace).

For the \textit{arietta} \textit{Ecco pur} in Act II, Perritt identifies strong cadential emphasis on the tonic even amid modulation. She notes that Monteverdi ``follows the correct rules for modulation," moving from $g$ minor to its relative major $B\flat$ and back to $g$ minor. The modulation occurs in measure 7 with a deceptive cadence in $B\flat$ major ($V$ moving to $vi$ instead of $I$), where the $vi$ chord is simultaneously the tonic of $g$ minor. She concludes that ``Monteverdi uses modulation to demonstrate tonality in this brief work.''

Our analysis reveals the deeper structural mechanism: the modulation between major and minor is explained by the transformation under $R$ (the relative operation) of the cadence $\{IV, V\}$ to $\{ii, iii\}$, which is preserved under this transformation as shown in our commutative diagrams. The dramatic function is particularly poignant: the section about the sun's return is entirely in major, while the actual \textit{ritorno} (return) is bittersweet—represented by a minor chord. Orpheus returns to the place that makes him happy, yet that same place harbors the sadness of Euridice's death.

Perritt's analysis of \textit{Mira, deh mira Orfeo} (also in Act II) notes that as the piece's mood changes, so does the tonality: joyful sections appear in $C$ major or closely related keys, while solemn sections shift to a minor (the relative minor of $C$ major). She emphasizes that ``the tonal system remains strong through the resolution of harmonic progressions,'' with secondary dominants ($V$/$V$) resolving properly to $V$ and then to $I$, calling this ``an interesting use of tonality that allows us to see the transition between modality and tonality.''

The modulation in this \textit{arietta} parallels that of \textit{Ecco pur}, alternating between major and minor tonalities. Sung by a shepherd, it foreshadows the drama to come: his happiness in major tonality is cut short by a messenger bringing the terrible news of Euridice's death. The shift to minor announces the catastrophe before it is explicitly stated.

Taken together, the modulations in these excerpts from \textit{L'Orfeo} narrate the drama Orpheus suffers for Euridice. This represents Monteverdi's practical response to Artusi's criticisms: the music serves the text and dramatic action, with tonal relationships deployed systematically to express emotional states and narrative progression. Where Perritt identifies the presence of tonal thinking in isolated passages, our quantum modulation model demonstrates how these modulations form a coherent dramatic-tonal architecture spanning entire scenes and acts.

\section{Conclusion and Future Work}
This article has demonstrated that algebraic models of counterpoint and modulation offer a powerful lens for analysing the structural genius of Claudio Monteverdi. By applying a symmetry-based framework to his counterpoint and group-theoretic models to his use of harmony, we have shown that his celebrated innovations were not arbitrary ``errors'' but part of a coherent, logical, and mathematically elegant expansion of musical language. The formalisms reveal a deep structural integrity behind the expressive force of the \textit{seconda pratica}, validating Monteverdi's artistic choices against the critiques of his time.
This work serves as a starting point for further inquiry. Future research could apply more advanced tools from algebra and group theory to explore the hidden mathematical structures in the broader repertoire of Monteverdi and his contemporaries. Such investigations promise to continue revealing the profound and often surprising connections between the aesthetic triumphs of music and the abstract beauty of mathematics.

%
%
%
\bibliographystyle{splncs04}
\bibliography{main}
\end{document}